\documentclass[11pt,a4paper]{amsart}
\usepackage{amsmath}
\usepackage{amssymb}
\usepackage{latexsym}
\usepackage{xypic}

\newcommand{\dual}{\makebox[0mm]{}^{{\scriptstyle\vee}}}
\newtheorem{theorem}{Theorem}[section]

\newtheorem{exmp}[theorem]{Example}
\newtheorem{exmps}[theorem]{Examples}
\newtheorem{rem}[theorem]{Remark}

\newcommand{\beeq}[1]{\begin{eqnarray}\label{#1}}
\newcommand{\eneq}{\end{eqnarray}}
\newcommand{\cal}{\mathcal}

\newcommand{\ka}{{\cal A}}

\newcommand{\kc}{{\cal C}}

\newcommand{\ke}{{\cal E}}

\newcommand{\kl}{{\cal L}}

\newcommand{\ko}{{\cal O}}

\newcommand{\Db}{{\rm D}^{\rm b}}

\newcommand{\IC}{{\mathbb C}}

\newcommand{\IP}{{\mathbb P}}
\newcommand{\IQ}{{\mathbb Q}}

\newcommand{\IZ}{{\mathbb Z}}

\newcommand{\Hom}{{\rm Hom}}
\newcommand{\Ext}{{\rm Ext}}

\newcommand{\verylongarrow}[1]{\hbox to #1{\rightarrowfill}}

\newcommand\mynote[1]{\marginpar{\ \\ \small \tt #1}}
\newcommand\bel[1]{{\mynote{#1}}\begin{equation}\label{#1}}

\begin{document}
\title{Proof of C\u{a}ld\u{a}raru's Conjecture.\\
An Appendix to a paper by K.\ Yoshioka}

\author[Daniel Huybrechts]{Daniel Huybrechts}
\address{D.H.: Mathematisches Institut, Universit{\"a}t Bonn,
  Beringstr.\ 1, D- 53115 Bonn, Germany}
\email{huybrech@math.uni-bonn.de}

\author[Paolo Stellari]{Paolo Stellari}
\address{P.S.: Dipartimento di Matematica ``F.
Enriques'', Universit{\`a} degli Studi di Milano, Via Cesare Saldini
50, 20133 Milano, Italy} \email{stellari@mat.unimi.it}

\maketitle
\pagestyle{plain}
In this short note we show how to combine Yoshioka's recent
results on moduli spaces of twisted sheaves on K3 surfaces with
more or less standard methods to prove C\u{a}ld\u{a}raru's
conjecture on the equivalence of twisted derived categories of
projective K3 surfaces. More precisely, we shall show

\begin{theorem} Let $X$ and $X'$ be two projective K3 surfaces
endowed with B-fields $B\in H^2(X,\IQ)$ respectively $B'\in
H^2(X',\IQ)$. Suppose there exists a Hodge isometry
$$g:\widetilde H(X,B,\IZ)\cong \widetilde H(X',B',\IZ)$$
that preserves the natural orientation of the four positive
directions. Then there exists a Fourier--Mukai equivalence
$\Phi:\Db(X,\alpha)\cong\Db(X',\alpha')$ such that the induced
action $\Phi^{B,B'}_*$ on cohomology equals $g$. Here,
$\alpha:=\alpha_B$ and $\alpha':=\alpha_{B'}$ are the Brauer classes
induced by $B$ respectively $B'$.
\end{theorem}

 The twisted Hodge structures and the cohomological Fourier-Mukai transform (based on the notion of twisted Chern character),
indispensable for the formulation of the conjecture, were introduced
in \cite{HS}. For a complete discussion of the natural orientation
of the positive directions and the cohomological Fourier--Mukai
transform $\Phi^{B,B'}_*$ we also refer to \cite{HS}. Note that
C\u{a}ld\u{a}raru's conjecture was originally formulated purely in
terms of the transcendental lattice. But, as has been explained in
\cite{HS}, in the twisted case passing from the transcendental part
to the full cohomology is not always possible, so that the original
formulation had to be changed slightly to the above one.

Also note that any Fourier--Mukai equivalence
$\Phi:\Db(X,\alpha)\cong\Db(X',\alpha')$ induces a Hodge isometry as
above, but for the time being we cannot prove that this Hodge
isometry also preserves the natural orientation. In the untwisted
case this is harmless, for a given orientation reversing Hodge
isometry can always be turned into an orientation preserving one by
composing with $-{\rm id}_{H^2}$. In the twisted setting this cannot
always be guaranteed, so that we cannot yet exclude the case of
Fourier--Mukai equivalent twisted K3 surfaces $(X,\alpha_B)$ and
$(X',\alpha_{B'})$ which only admit an orientation reversing Hodge
isometry $\widetilde H(X,B,\IZ)\cong \widetilde H(X',B',\IZ)$. Of
course, this is related to the question whether any Fourier--Mukai
equivalence is orientation preserving which seems to be a difficult
question even in the untwisted case (see \cite{HLOY,Sz}).

From Yoshioka's paper \cite{Y} we shall use the following

\begin{theorem}{\bf (Yoshioka)}\label{Yoshthm}
Let $X$ be a K3 surface with $B\in H^2(X,\IQ)$ and $v\in
\widetilde H^{1,1}(X,B,\IZ)$ a primitive vector with $\langle
v,v\rangle=0$. Then there exists a moduli space $M(v)$ of stable
(with respect to a generic polarizations) $\alpha_B$-twisted
sheaves $E$ with ${\rm ch}^B(E) \sqrt{td(X)}=v$  such that:

{\rm i)} Either $M(v)$ is empty or a K3 surface. The latter holds
true if the degree zero part of $v$ is positive.

{\rm ii)} On $X':=M(v)$ one finds a B-field $B'\in H^2(X',\IQ)$
such that there exists a universal family $\ke$ on $X\times X'$
which is an $\alpha_B^{-1}\boxtimes\alpha_{B'}$-twisted sheaf.

{\rm iii)} The twisted sheaf $\ke$ induces a Fourier--Mukai
equivalence $\Db(X,\alpha_B)\cong\Db(X',\alpha_{B'})$.
\end{theorem}

The existence of the moduli space of semistable twisted sheaves  has
been proved by Yoshioka for arbitrary projective varieties. Instead of
considering twisted sheaves, he works with coherent sheaves on a
Brauer--Severi variety. Using the equivalence between twisted sheaves
and modules over Azumaya algebras, one can in fact view these moduli
spaces also as a special case of Simpson's general construction
\cite{Simpson}. (The two stability conditions are indeed equivalent.)
In his thesis \cite{L} M.\ Lieblich considers similar moduli spaces. (See also
\cite{HofSt} for the rank one case.)

The crucial part for the application to  C\u{a}ld\u{a}raru's
 conjecture is i), in particular the non-emptiness. Yoshioka follows Mukai's approach, which also
 yields ii). Part iii) is a rather formal consequence of the usual criteria for the
equivalence of Fourier--Mukai transforms already applied to the
twisted case in \cite{Cal}.

In the last section we provide a dictionary between the
different versions of twisted Chern characters and the various
notions of twisted sheaves. Only parts of it is actually used in
the proof of the conjecture. The rest is meant to complement
\cite{HS} and to facilitate the comparison of \cite{HS}, \cite{L},
and \cite{Y}.
\medskip

{\bf Acknowledgements:} It should be clear that the lion's share of
the proof of C\u{a}ld\u{a}raru's conjecture in the above form is in
fact contained in K.\ Yoshioka's paper. We are grateful to him for
informing us about his work and comments on the various versions of
this note. Thanks also to M.\ Lieblich who elucidated the relation
between the different ways of constructing moduli spaces of twisted
sheaves. Our proof follows the arguments in the untwisted case, due
to S.\ Mukai \cite{Mukai} and D.\ Orlov \cite{Orlov} (see also
\cite{HLOY,Pl}), although some modifications were necessary. During
the preparation of this paper the second named author was partially
supported by the MIUR of the Italian government in the framework of
the National Research Project ``Algebraic Varieties'' (Cofin 2002).

\section{Examples}\label{Exas}

Let $X$ and $X'$ be projective K3 surfaces (always over $\IC$) with
B-fields $B\in H^2(X,\IQ)$ respectively $B'\in H^2(X',\IQ)$. We
denote the induced Brauer classes by
$\alpha:=\alpha_B:=\exp(B^{0,2})\in H^2(X,\ko_X^*)$ respectively
$\alpha':=\alpha_{B'}\in H^2(X',\ko_{X'}^*)$. We start out with
introducing a few examples of equivalences between the bounded
derived categories $\Db(X,\alpha)$ respectively $\Db(X',\alpha')$ of
the abelian categories of $\alpha$-twisted (resp.\
$\alpha'$-twisted) sheaves.

\medskip

i) Let $f:X\cong X'$ be an automorphism with $f^*\alpha'=\alpha$.
Then $\Phi:=f_*:\Db(X,\alpha)\to\Db(X',\alpha'), E\mapsto Rf_*E$ is
a Fourier--Mukai equivalence with kernel $\ko_{\Gamma_f}$ viewed as
an $\alpha^{-1}\boxtimes\alpha'$-twisted sheaf on $X\times X'$.

If in addition $f_*(B)=B'$ then $\Phi^{B,B'}_*=f_*$.

\medskip

ii) Let $L\in{\rm Pic}(X)$ be a(n untwisted) line bundle on $X$.
Then $E\mapsto L\otimes E$ defines a Fourier--Mukai equivalence
$L\otimes(~~):\Db(X,\alpha)\cong\Db(X,\alpha)$ with kernel $i_*L$
considered as an $\alpha^{-1}\boxtimes\alpha$-twisted sheaf on
$X\times X$. Here, $i:X\hookrightarrow X\times X$ denotes the
diagonal embedding. The induced cohomological Fourier--Mukai
transform $(L\otimes(~~))^{B,B}_*:\widetilde H(X,B,\IZ)\cong
\widetilde H(X,B,\IZ)$ is given by multiplication with $\exp({\rm
c}_1(L))$.

\medskip

iii) Let $b\in H^2(X,\IZ)$. Then $\alpha_B=\alpha_{B+b}$. The
identity ${\rm id}:\Db(X,\alpha_B)=\Db(X,\alpha_{B+b})$ descends to
${\rm id}^{B,B+b}_*:\widetilde H(X,B,\IZ)\cong \widetilde
H(X,B+b,\IZ)$ which is given by the multiplication with $\exp(b)$.
This follows from the formula ${\rm ch}^{B+b}(E)={\rm
ch}^B(E)\cdot{\rm ch}^b(\ko)={\rm ch}^B(E)\cdot\exp(b)$ (see
\cite[Prop.\ 1.2]{HS}).

\medskip

iv) Changing the given B-field $B$  by a class $b\in
H^{1,1}(X,\IQ)$ does not affect $\widetilde H(X,B,\IZ)$. Thus,
the identity can be considered as an orientation preserving Hodge
isometry $\widetilde H(X,B,\IZ)=\widetilde H(X,B+b,\IZ)$.

As shall be explained in the last section, this can be lifted to a
Fourier--Mukai equivalence. More precisely, there is an exact
functor $\Phi :{\rm\bf Coh}(X,\alpha_B)\cong{\rm\bf
Coh}(X,\alpha_{B+b})$, whose derived functor, again denoted by
$\Phi:\Db(X,\alpha_B)\cong\Db(X,\alpha_{B+b})$, is of Fourier--Mukai
type and such that $\Phi^{B,B+b}_*={\rm id}$.

\medskip

v) Let $E\in\Db(X,\alpha)$ be a spherical object, i.e.\
$\Ext^i(E,E)=0$ for all $i$ except for $i=0,2$ when it is of
dimension one. Then the twist functor $T_E$ that sends
$F\in\Db(X,\alpha)$ to the cone of $\Hom(E,F)\otimes E\to F$ defines
a Fourier--Mukai autoequivalence
$T_E:\Db(X,\alpha)\cong\Db(X,\alpha)$. The kernel of $T_E$ is given
by the cone of the natural map
$$E^*\boxtimes E\to\ko_\Delta,$$
where $\ko_\Delta$ is considered as an
$\alpha^{-1}\boxtimes\alpha$-twisted sheaf on $X\times X$. The
result in the untwisted case goes back to Seidel and Thomas
\cite{ST}. The following short proof of this, which carries over to
the twisted case, has been communicated to us by D.\ Ploog
\cite{Pl}. Consider the class $\Omega\subset\Db(X,\alpha)$ of
objects $F$ that are either isomorphic to $E$ or contained in its
orthogonal complement $E^\perp$, i.e.\ $\Ext^i(E,F)=0$ for all $i$.
It is straightforward to check that this class is spanning. Since
$T_E(E)\cong E[-1]$ and $T_E(F)\cong F$ for $F\in E^\perp$, one
easily verifies that $\Ext^i(F_1,F_2)=\Ext^i(T_E(F_1),T_E(F_2))$ for
all $F_1,F_2\in\Omega$.

In other words, $T_E$ is fully faithful on the spanning class
$\Omega$ and hence fully faithful. By the usual argument, the
Fourier--Mukai functor $T_E$ is then an equivalence.

As in the untwisted case, one proves that the induced action on
cohomology is the reflection $\alpha\mapsto\alpha
+\langle\alpha,v^B(E)\rangle\cdot v^B(E)$. Here, $v^B(E)$ is the
Mukai vector $v^B(E):={\rm ch}^B(E)\sqrt{{\rm td}(X)}$.

Special cases of this construction are:

 -- Let
$\IP^1\cong C\subset X$ be a smooth rational curve. As
$H^2(C,\ko_C^*)$ is trivial, its structure sheaf $\ko_C$ and any
twist $\ko_C(k)$ can naturally be considered as $\alpha$-twisted
sheaves. The Mukai vector for $k=-1$ is given by
$v(\ko_C(-1))=(0,[C],0)$.

-- In the untwisted case, the trivial line bundle $\ko$ (and in fact
any line bundle) provides an example of a spherical object. Its
Mukai vector is $(1,0,1)$ and has, in particular, a non-trivial
degree zero component. It is the latter property that is of
importance for the proof in the untwisted case. So the original
argument goes through if at least one spherical object of
non-trivial rank can be found. Unfortunately, spherical object (in
particular those of positive rank) might not exist at all in the
twisted case. In fact, any spherical object $E$ has a Mukai vector
$v^B(E)\in{\rm Pic}(X,B)$ of square $\langle
v^B(E),v^B(E)\rangle=-2$ and it is not difficult to find examples of
rational B-fields $B\ne0$ such that such a vector does not exist.

\medskip

vi) Let $\ell\in{\rm Pic}(X)$ be a nef class with
$\langle\ell,\ell\rangle=0$. If $w=(0,\ell,s)$ is a primitive
vector, then the moduli space $M(w)$ of $\alpha_B$-twisted sheaves
which are  stable with respect to a generic polarization is
non-empty. Indeed, in this case $\ell$ is a multiple $n\cdot f$ of a
fibre class $f$ of an elliptic fibration $\pi:X\to\IP^1$. As
$\gcd(n,s)=1$, there exists a stable rank $n$ vector bundle of
degree $s$ on a smooth fibre of $\pi$ which yields a point in
$M(w)$.

If $\ell$ is the fibre class of an elliptic fibration $X\to\IP^1$,
we can think of $M(w)$ as the relative Jacobian ${\cal
  J}^s(X/\IP^1)\to\IP^1$.

In any case, $M(w)$ is a K3 surface and the universal twisted sheaf
provides an equivalence $\Phi:\Db(M(w),\alpha_{B'})\cong
\Db(X,\alpha_B)$ (for some B-field $B'$ on $M(w)$) inducing a Hodge
isometry $\Phi^{B',B}_*:\widetilde H(M(w),B',\IZ)\cong \widetilde
H(X,B,\IZ)$ that sends $(0,0,1)$ to $w$.

\section{The proof}

Let $g:\widetilde H(X,B,\IZ)\cong \widetilde H(X',B',\IZ)$ be an
orientation preserving Hodge isometry. The Mukai vector of $k(x)$
with $x\in X$ is $v^B(k(x))=v(k(x))=(0,0,1)$. We shall denote its
image under $g$ by $w:=g(0,0,1)=(r,\ell,s)$.

\medskip

{\bf 1st step.} In the first step we assume that $r=0$ and
$\ell=0$, i.e.\ $g(0,0,1)=\pm(0,0,1)$, and that furthermore
$g(1,0,0)=\pm(1,0,0)$. By composing with $-{\rm id}$ we may actually
assume $g(0,0,1)=(0,0,1)$ and $g(1,0,0)=(1,0,0)$.

In particular, $g$ preserves the grading of $\widetilde H$ and
induces a Hodge isometry $H^2(X,\IZ)\cong H^2(X',\IZ)$. Denote
$b:=g(B)-B'\in H^2(X,\IQ)$. As $g$ respects the Hodge structure, it
maps $\sigma+B\wedge \sigma$ to $\sigma'+B'\wedge\sigma'$ and,
therefore $\langle\sigma,B\rangle=\langle\sigma',B'\rangle$. On the
other hand, as $g$ is an isometry, one has
$\langle\sigma,B\rangle=\langle\sigma',g(B)\rangle$. Altogether this
yields $\langle\sigma',b\rangle=0$, i.e.\ $b\in H^{1,1}(X,\IQ)$.

Now compose $g$ with the orientation preserving Hodge isometry given
by the identity $\widetilde H(X',B',\IZ)=\widetilde
H(X',g(B)=B'+b,\IZ)$. As the latter can be lifted to a
Fourier--Mukai equivalence $\Db(X',\alpha_{B'})\cong
\Db(X',\alpha_{g(B)})$ (see example iv)), it suffices to show that
$g$ viewed as a Hodge isometry $\widetilde H(X,B,\IZ)=\widetilde
H(X,g(B),\IZ)$ can be lifted. So we may from now on assume that
$B'=g(B)$.

As $g$ is orientation preserving, its degree two component defines a
Hodge isometry that maps the positive cone $\kc_X\subset H^{1,1}(X)$
onto the positive cone $\kc_{X'}\subset H^{1,1}(X')$.

If $g$ maps an ample class to an ample class, then by the Global
Torelli Theorem $g$ can be lifted to an isomorphism $f:X\cong X'$
which in turn yields  a Fourier--Mukai equivalence
$\Phi:=f_*:\Db(X,\alpha_B)\cong\Db(X',\alpha_{B'})$. Obviously,
with this definition $\Phi^{B,B'}_*=g$ (use $f_*(B)=g(B)=B'$).

If $g$ does not preserve ampleness, then the argument has to be
modified as follows: After a finite number of reflections
$s_{C_i}$ in hyperplanes orthogonal to $(-2)$-classes $[C_i]$ we
may assume that $s_{C_1}(\ldots s_{C_n}(h(a))\ldots)$ is an ample
class. As the reflections $s_{C_i}$ are induced by the twist
functors
$T_{\ko_{C_i}(-1)}:\Db(X',\alpha_{B'})\cong\Db(X',\alpha_{B'})$
(see the explanations in the last section), the Hodge isometry
$g$ is induced by a Fourier--Mukai equivalence if and only if the
composition $s_{C_1}\circ \ldots s_{C_n}\circ g$ is. Thus, we
have reduced the problem to the case already treated above.

In the following steps we shall explain how the general case can
be reduced to the case just considered.

\medskip

{\bf 2nd step.} Suppose $g(0,0,1)=\pm(0,0,1)$ but
$g(1,0,0)\ne\pm(1,0,0)$.
Again, by composing with $-{\rm id}$ we may reduce to
$g(0,0,1)=(0,0,1)$
and $g(1,0,0)\ne(1,0,0)$.
Then $g(1,0,0)$ is necessarily of the form $\exp(b)$ for some
$b\in H^2(X',\IZ)$. Hence, we may compose $g$ with the  Hodge
isometry $\exp(-b):\widetilde H(X',B',\IZ)\cong \widetilde
H(X',B'-b,\IZ)$ (that preserves the orientation) which can be
lifted to a Fourier--Mukai equivalence according to example iii).
This reduces the problem to the situation studied in the previous
step.

\medskip

{\bf 3rd step.}
 Suppose that $r>0$. Using Theorem \ref{Yoshthm}
one finds a K3 surface $X_0$ with a B-field $B_0\in H^2(X_0,\IQ)$
such that over $X_0\times X'$ there exists a universal
$\alpha_{B_0}^{-1}\boxtimes\alpha_{B'}$-twisted sheaf parametrizing
stable $\alpha'$-twisted sheaves on $X'$ with Mukai vector
$v^{B'}=w$. In particular, $\ke$ induces an equivalence
$\Phi_\ke:\Db(X_0,\alpha_{B_0})\cong\Db(X',\alpha_{B'})$ and
$\Phi^{B_0,B'}_{\ke*}(0,0,1)=w$.

Thus, the composition $g_0:=(\Phi^{B_0,B'}_{\ke*})^{-1}\circ g$
yields an orientation preserving (!) Hodge isometry $\widetilde
H(X,B,\IZ)\cong H(X_0,B_0,\IZ)$. (The proof that the universal
family of stable sheaves induces an orientation preserving Hodge
isometry is analogous to the untwisted case. This seems to be widely
known \cite{HLOY,Sz}. For an explicit proof see \cite{HS}.) Clearly,
$g$ can be lifted to a Fourier--Mukai equivalence if and only if
$g_0$ can. The latter follows from step one.

\medskip

{\bf 4th step.} Suppose $g$ is given with $r<0$. Then compose with
the orientation preserving Hodge isometry $-{\rm id}$ of
$\widetilde H(X',B',\IZ)$ which is lifted to the shift functor
$E\mapsto E[1]$. Thus, it is enough to lift the composition $-{\rm
id}\circ g$ which can be achieved according to step three.

\bigskip

{\bf 5th step} The remaining case is $r=0$ and $\ell\ne0$. One
applies the construction of example vi) in Section \ref{Exas} and
proceeds as in step 3. The class $\ell$ can be made nef by applying
$-{\rm id}$ if necessary to make it effective (i.e.\ contained in
the closure of the positive cone) and then composing it with
reflections $s_C$ as in step one.

\section{The various twisted categories and their Chern characters}
Let $\alpha\in H^2(X,\ko_X^*)$ be a Brauer class represented by
a \v{C}ech cocycle $\{\alpha_{ijk}\}$.

{\bf 1.} The abelian category ${\rm\bf Coh}(X,\{\alpha_{ijk}\})$
of $\{\alpha_{ijk}\}$-twisted coherent sheaves
only depends on the class $\alpha\in H^2(X,\ko_X^*)$. More
precisely, for any other choice of a \v{C}ech-cocycle
$\{\alpha'_{ijk}\}$ representing $\alpha$, there exists an
equivalence
$$\xymatrix{{\rm\bf
Coh}(X,\{\alpha_{ijk}\})\ar[r]^-{\Psi_{\{\lambda_{ij}\}}}&{\rm\bf
Coh}(X,\{\alpha'_{ijk}\}),~~
\{E_i,\varphi_{ij}\}\ar@{|->}[r]&\{E_i,\varphi_{ij}\cdot\lambda_{ij}\},}$$
where $\{\lambda_{ij}\in\ko^*(U_{ij})\}$ satisfies
$\alpha'_{ijk}\alpha_{ijk}^{-1}=\lambda_{ij}\cdot\lambda_{jk}\cdot\lambda_{ki}$.
Clearly, $\{\lambda_{ij}\}$ exists, as $\{\alpha_{ijk}\}$ and
$\{\alpha'_{ijk}\}$ define the same Brauer class, but it is far from
being unique. In other words, the above equivalence
$\Psi_{\{\lambda_{ij}\}}$ is not canonical. In order to make this
more precise, choose a second $\{\lambda'_{ij}\}$. Then
$\gamma_{ij}:=\lambda_{ij}'\cdot\lambda_{ij}^{-1}$ can be viewed as
the transition function of a holomorphic line bundle
$\kl_{\lambda\lambda'}$. With this notation one finds
$$\Psi_{\{\lambda_{ij}'\}}=(\kl_{\lambda\lambda'}\otimes(\;\;))
\circ\Psi_{\{\lambda_{ij}\}}.$$

A very special case of this is the equivalence
$$\xymatrix{\kl\otimes(\;\;):{\rm\bf
Coh}(X,\{\alpha_{ijk}\})\ar[r]&{\rm\bf Coh}(X,\{\alpha_{ijk}\})}$$
that is induced by the tensor product with a holomorphic line bundle
$\kl$ given by a cocycle $\{\gamma_{ij}\}$.

Despite this ambiguity in identifying these categories for
different choices of the \v{C}ech-representative, ${\rm\bf
Coh}(X,\{\alpha_{ijk}\})$ is often simply denoted ${\rm\bf
Coh}(X,\alpha)$.

\medskip

{\bf  2.} Now fix a B-field $B\in H^2(X,\IQ)$ together with a
\v{C}ech-representative $\{B_{ijk}\}$. The induced Brauer class
$\alpha:=\exp(B^{0,2})\in H^2(X,\ko_X^*)$ is represented by the
\v{C}ech-cocycle $\{\alpha_{ijk}:=\exp(B_{ijk})\}$.

In \cite{HS} we introduced
$$\xymatrix{{\rm ch}^B:{\rm \bf
Coh}(X,\{\alpha_{ijk}\})\ar[r]&H^{*,*}(X,\IQ).}$$

The construction makes use of a further choice of ${\cal
C}^\infty$-functions $a_{ij}$ with
$-B_{ijk}=a_{ij}+a_{jk}+a_{ki}$, but the result does not depend
on it. Indeed, by definition, ${\rm
ch}^B(\{E_i,\varphi_{ij}\})={\rm
ch}(\{E_i,\varphi_{ij}\cdot\exp(a_{ij})\})$. Thus, if we pass
from $a_{ij}$ to $a_{ij}+c_{ij}$ with $c_{ij}+c_{jk}+c_{ki}=0$,
then ${\rm ch}^B(\{E_{ij},\varphi_{ij}\})$ changes by $\exp({\rm
c}_1(\kl))$, where $\kl$ is given by the transition functions
$\{\exp(c_{ij})\}$. But by the very definition of the first Chern
class, one has ${\rm c}_1(\kl)=[\{c_{ij}+ c_{jk}+ c_{ki}\}]=0$.

More generally, we may change the class $B$ by a class $b\in
H^2(X,\IQ)$ represented by $\{b_{ijk}\}$. Suppose
$\alpha_{B+b}=\alpha_{B}\in H^2(X,\ko_X^*)$. We denote the
\v{C}ech-representative $\exp(B_{ijk}+b_{ijk})$ by
$\alpha'_{ijk}$. As before, we write
$-B_{ijk}=a_{ij}+a_{jk}+a_{ki}$ and
$-b_{ijk}=c_{ij}+c_{jk}+c_{ki}$.

The Chern characters ${\rm ch}^B$ and ${\rm ch}^{B+b}$ fit into the
following commutative diagram
$$\xymatrix{{\rm\bf
Sh}(X,\{\alpha_{ijk}\})\ar[rr]^{\Psi_{\{\exp(-c_{ij})\}}}\ar[dr]_{{\rm
ch}^B}&&{\rm\bf Sh}(X,\{\alpha'_{ijk}\})\ar[dl]^{{\rm
ch}^{B+b}}\\
&H^*(X,\IQ).&}$$ Unfortunately, we cannot replace ${\rm\bf Sh}$ by
${\rm\bf Coh}$, for $\exp(c_{ij})$ are only differen\-tiable
functions. Nevertheless, there exist non-unique
$\beta_{ij}\in\ko^*(U_{ij})$ with
$\alpha_{ijk}'=\alpha_{ijk}\cdot(\beta_{ij}\cdot\beta_{jk}\cdot\beta_{ki})$.
Using these one finds a commutative diagram
$$\xymatrix{{\rm\bf
Coh}(X,\{\alpha_{ijk}\})\ar[d]_{{\rm ch}^B}
\ar[rr]^{\Psi_{\{\beta_{ij}\}}}&&
{\rm\bf Coh}(X,\{\alpha'_{ijk}\})\ar[d]^{{\rm ch}^{B+b}}\\
H^*(X,\IQ)\ar[rr]_{\exp({\rm c}_1(\kl))}&&H^*(X,\IQ). }$$ Here,
${\kl}$ is the line bundle given by the transition functions
$\beta_{ij}\cdot\exp(c_{ij})$.

It is not difficult to see that ${\rm c}_1(\kl)\in H^{1,1}(X)$
whenever $b\in H^{1,1}(X,\IQ)$. Indeed, ${\rm
c}_1(\kl)=\{d\log(\beta_{ij})\}+b$, which is of type $(1,1)$, as
$b$ is $(1,1)$ by assumption and the functions $\beta_{ij}$ are
holomorphic.

Thus, in this case there exists a holomorphic line bundle
$\widetilde\kl$ with ${\rm c}_1(\widetilde \kl)={\rm c}_1(\kl)$. Now
consider the composition
$\Phi:=(\widetilde\kl^*\otimes(~~))\circ\Psi_{\{\beta_{ij}\}}:
{\rm\bf Coh}(X,\alpha_B:=\{\alpha_{ijk}\})\cong{\rm\bf
Coh}(X,\alpha_{B+b}:=\{\alpha'_{ijk}\})$, which is an exact
equivalence, and denote the derived one again by
$\Phi:\Db(X,\alpha_B)\cong\Db(X,\alpha_{B+b})$. Then the above
calculation of the twisted Chern character implies that
$\Phi^{B,B+b}_*={\rm id}$.

\medskip

{\bf 3.} Consider again the abelian category ${\rm\bf
Coh}(X,\{\alpha_{ijk}\})$. For any locally free
$G=\{G_i,\varphi_{ij}\}\in{\rm\bf Coh}(X,\{\alpha_{ijk}\})$ one
defines an Azumaya algebra $\ka_G:=\ke nd (G\dual)$. The abelian
category of left $\ka_G$-modules will be denoted ${\rm\bf
Coh}(\ka_G)$. An equivalence of abelian categories is given by
$$\xymatrix{{\rm\bf Coh}(X,\{\alpha_{ijk}\})\ar[r]&{\rm\bf Coh}(\ka_G),~~
E\ar@{|->}[r]&G\dual\otimes E.}$$

In \cite{Y} Yoshioka considers yet another abelian category
${\rm\bf Coh}(X,Y)$ of certain coherent sheaves on a projective
bundle $Y\to X$ realizing the Brauer class $\alpha$. As is
explained in detail in \cite{Y}, one again has an equivalence of
abelian categories ${\rm\bf Coh}(X,Y)\cong{\rm\bf
Coh}(X,\{\alpha_{ijk}\})$. In order to define an appropriate
notion of stability, Yoshioka defines a Hilbert polynomial  for
objects $E\in{\rm\bf Coh}(X,Y)$. It is straightforward to see that
under the composition
$$\xymatrix{{\rm\bf Coh}(X,Y)\ar[r]^-\sim&{\rm\bf
Coh}(X,\{\alpha_{ijk}\})\ar[r]^-\sim&{\rm \bf Coh}(\ka_G)}$$ his
Hilbert polynomial corresponds to the usual Hilbert polynomial
for sheaves $F\in{\rm\bf Coh}(\ka_G)$ viewed as $\ko_X$-modules.
The additional choice of the locally free object $G$ in ${\rm\bf
Coh}(X,Y)$ or equivalently in ${\rm\bf Coh}(X,\{\alpha_{ijk}\})$
needed to define the Hilbert polynomial in \cite{Y} enters this
comparison via the equivalence ${\rm\bf
Coh}(X,\{\alpha_{ijk}\})\cong{\rm \bf Coh}(\ka_G)$. From here it
is easy to see that the stability conditions considered in
\cite{Simpson,Y}  are actually equivalent.

 We would like to define a twisted
Chern character for objects in ${\rm\bf Coh}(\ka_G)$. Of course,
as any $F\in{\rm\bf Coh}(\ka_G)$ is in particular an ordinary
sheaf, ${\rm ch}(F)$ is well defined. In order to define
something that takes into account the $\ka_G$-module structure,
one fixes $B=\{B_{ijk}\}$ and assumes
$\alpha_{ijk}=\exp(B_{ijk})$. Then we introduce
$$\xymatrix{{\rm ch}^B_G:{\rm\bf
Coh}(\ka_G)\ar[r]&H^*(X,\IQ),~~F\ar@{|->}[r]&\frac{{\rm ch}(F)}{{\rm
ch}^{-B}(G\dual)}.}$$ Note that {\it a priori} the definition
depends on $B$ and $G$, but the dependence on the latter is
well-behaved as will be explained shortly.

Here are the main compatibilities for this new Chern character:

{\bf i)} The following diagram is commutative:
$$\xymatrix{{\rm\bf Coh}(X,\{\alpha_{ijk}\})\ar[dr]_{{\rm
ch}^B}\ar[rr]&&
{\rm\bf Coh}(\ka_G)\ar[dl]^{{\rm ch}_G^B}\\
&H^*(X,\IQ).&}$$ Indeed, ${\rm ch}(G\dual\otimes E)={\rm
ch}^{-B}(G\dual)\cdot{\rm ch}^B(E)$.

{\bf ii)} Let $H$ be a locally free coherent sheaf and
$G':=G\otimes H\in{\rm\bf Coh}(X,\{\alpha_{ijk}\})$.  Then the
natural equivalence ${\rm\bf Coh}(\ka_G)\to{\rm\bf
Coh}(\ka_{G'})$, $F\mapsto H\dual\otimes F$ fits in the
commutative diagram
$$\xymatrix{{\rm\bf Coh}(\ka_G)\ar[dr]_{{\rm ch}_G^B}
\ar[rr]&&{\rm\bf Coh}(\ka_{G'})\ar[dl]^{{\rm ch}^B_{G'}}\\
&H^*(X,\IQ).&}$$ This roughly says that the new Chern character is
independent of $G$.

{\bf iii)} If $E_1,E_2\in {\rm\bf Coh}(X,\{\alpha_{ijk}\})$ and
$F_i:=G\dual\otimes E_i\in{\rm\bf Coh}(\ka_G)$ then
$\chi(E_1,E_2):=\sum(-1)^i\dim{\rm Ext}^i(E_1,E_2)$ is well-defined
and equals $\chi(F_1,F_2):=\sum(-1)^i\dim{\rm
Ext}^i_{\ka_G}(F_1,F_2)$. Both expressions can be computed in terms
of the twisted Chern characters introduced above and the Mukai
pairing. Concretely,
$$\chi(F_1,F_2)=-\langle{\rm ch}_G^B(F_1)\cdot\sqrt{{\rm td}(X)},
{\rm ch}_G^B(F_2)\cdot\sqrt{{\rm td}(X)}\rangle.$$ Here
$\langle~~,~~\rangle$ denotes the generalized Mukai pairing and
$$\chi(F_1,F_2):=-\langle{\rm ch}^B(E_1)\cdot\sqrt{{\rm td}(X)},{\rm ch}^B(E_2)\cdot\sqrt{{\rm td}(X)}\rangle.$$ (Be aware of the different sign conventions for
K3 surfaces and the general case.)

\medskip

{\bf 4.} There is yet another way to define a twisted Chern
character which is implicitly used in \cite{Y}. We use the above
notations and define ${\rm ch}_G:{\rm \bf Coh}(\ka_G)\to
H^*(X,\IQ)$ by ${\rm ch}_G(F):=\frac{{\rm ch}(F)}{\sqrt{{\rm
ch}(\ka_G)}}$, where $F$ and $\ka_G$ are considered as ordinary
$\ko_X$-modules. Using the natural identifications explained
earlier, namely ${\rm \bf Coh}(X,Y)\cong{\rm\bf
Coh}(X,\{\alpha_{ijk}\})\cong{\rm\bf Coh}(\ka_G)$, this Chern
character can also be viewed as a Chern character on the other
abelian categories.

Although the definition ${\rm ch}_G$ seems very natural, it does
not behave nicely under change of $G$. More precisely, in
general  ${\rm ch}_{G\otimes H}(H\dual\otimes F)\ne{\rm ch}_G(F)$.

Fortunately, the situation is less critical for K3 surfaces.
Here, the relation between ${\rm ch}_G$ and ${\rm ch}_G^B$ can be
described explicitly and using the results in {\bf 3.} one deduces
from this a formula for the change of ${\rm ch}_G$ under
$G\mapsto G\otimes H$. In fact, it is straightforward to see that
the following diagram commutes:
$$\xymatrix{&{\rm \bf Coh}(\ka_G)\ar[dl]_{{\rm
ch}_G^B}\ar[dr]^{{\rm ch}_G}&\\
H^*(X,\IQ)\ar[rr]_{\exp(-B_G)}&&H^*(X,\IQ).}$$ Here $B_G:=\frac{{\rm
c}_1^B(G)}{{\rm rk}(G)}$, where ${\rm c}_1^B(G)$ is the degree two
part of ${\rm ch}^B(G)$. Note that $B$ and $B_G$ define the same
Brauer class. In particular, the Hodge structures $\widetilde
H^*(X,B,\IZ)$ and $\widetilde H^*(X,B_G,\IZ)$ are isomorphic.

This relation between ${\rm ch}_G^B$ and ${\rm ch}_G$ can be used
to compare the two versions of the cohomogical Fourier--Mukai
transform in \cite{HS} and \cite{Y}. With $v^B:={\rm
ch}^B\cdot\sqrt{{\rm td}(X)}$ and $v_G:={\rm ch}_G\cdot\sqrt{{\rm
td}(X)}$ and the implicit identification ${\rm \bf
Coh}(X,\alpha)={\rm \bf Coh}(\ka_G)$ the following diagram is
commutative:
$$\xymatrix{\Db(X,\alpha_B)\ar[rrr]^{\Phi}\ar[dd]_{v^B}\ar[dr]^{v_G}&&&
\Db(X',\alpha')\ar[dd]^{v^{B'}}\ar[dl]_{v_{G'}}\\
&H^*(X,\IQ)\ar[r]&H^*(X',\IQ)&\\
H^*(X,B,\IZ)\ar[ur]_{\exp(-B_G)}\ar[rrr]_{\Phi_*^{B,B'}}&&&H^*(X',B',\IZ).\ar[ul]^{\exp(-B_{G'})}}$$
Here, the central isomorphism $H^*(X,\IQ)\cong H^*(X',\IQ)$ is the
correspondence defined  by $v_{G\dual\boxtimes G'}(\ke)$ with
$\ke\in\Db(X\times X',\alpha^{-1}_B\boxtimes\alpha_{B'})$ the kernel
defining $\Phi$.


\end{document}